\documentclass[12pt,twoside]{article}
\usepackage[english]{babel}
\usepackage[latin1]{inputenc}
\usepackage{amsmath}
\usepackage{amssymb,amsfonts}
\usepackage{graphicx}

\newcommand{\bG}{\mathbf{G}}

\newcommand{\bL}{\mathbf{L}}

\newcommand{\bR}{\mathbf{R}}

\newcommand{\bS}{\mathbf{S}}

\newcommand{\ba}{\mathbf{a}}
\newcommand{\bc}{\mathbf{c}}
\newcommand{\be}{\mathbf{e}}

\newcommand{\bl}{\mathbf{l}}
\newcommand{\bh}{\mathbf{h}}
\newcommand{\bj}{\mathbf{j}}
\newcommand{\bq}{\mathbf{q}}

\newcommand{\HYP}{\mathbb{H}^3}
\newcommand{\HYN}{\mathbb{H}^n}

\newcommand{\SLR}{\widetilde{\bS\bL_2\bR}}
\newcommand{\NIL}{\mathbf{Nil}}
\newcommand{\SOL}{\mathbf{Sol}}

\begin{document}
\pagestyle{myheadings}
\markboth{\centerline{Jen\H o Szirmai}}
{Congruent and non-congruent hyperball packings $\dots$}
\title
{Congruent and non-congruent hyperball packings related to doubly truncated Coxeter orthoschemes in hyperbolic $3$-space}

\author{\normalsize{Jen\H o Szirmai} \\
\normalsize Budapest University of Technology and \\
\normalsize Economics, Institute of Mathematics, \\
\normalsize Department of Geometry \\
\date{\normalsize{\today}}}

\maketitle


\begin{abstract}

In \cite{Sz17-2} we considered hyperball packings in $3$-dimensional hyperbolic space. We developed a decomposition algorithm that 
for each saturated hyperball packing provides a decomposition of $\HYP$ into truncated tetrahedra. 
In order to get a density upper bound for hyperball packings, it is sufficient to determine
the density upper bound of hyperball packings in truncated simplices. Therefore, in this paper we examine the doubly truncated Coxeter 
orthoscheme tilings and the corresponding congruent and non-congruent hyperball packings. We proved that related to the mentioned Coxeter tilings 
the density of the densest congruent hyperball packing is $\approx 0.81335$ that is -- by our conjecture -- the upper bound density of the relating non-congruent hyperball packings, too.
\end{abstract}

\newtheorem{theorem}{Theorem}[section]
\newtheorem{corollary}[theorem]{Corollary}
\newtheorem{conjecture}{Conjecture}[section]
\newtheorem{lemma}[theorem]{Lemma}
\newtheorem{exmple}[theorem]{Example}
\newtheorem{defn}[theorem]{Definition}
\newtheorem{rmrk}[theorem]{Remark}
\newenvironment{definition}{\begin{defn}\normalfont}{\end{defn}}
\newenvironment{remark}{\begin{rmrk}\normalfont}{\end{rmrk}}
\newenvironment{example}{\begin{exmple}\normalfont}{\end{exmple}}
\newenvironment{acknowledgement}{Acknowledgement}


\section{Introduction}
In $n$-dimensional hyperbolic space $\mathbb{H}^n$ $(n\ge2)$ there are $3$ kinds
of ''balls (spheres)": the classical balls (spheres), horoballs (horospheres) and hyperballs (hyperspheres).

In this paper we consider the hyperballs and their packings related to the complete Coxeter tilings in $3$-dimensional hyperbolic space $\HYP$. 
However, first we survey the previous results related to this topic.

In the hyperbolic plane $\mathbb{H}^2$ the universal upper bound of the hypercycle packing density is $\frac{3}{\pi}$,
proved by I.~Vermes in \cite{V79} and the universal lower bound of the hypercycle covering density is $\frac{\sqrt{12}}{\pi}$
determined by I.~Vermes in \cite{V81}. 

In \cite{Sz06-1} and \cite{Sz06-2} we analysed the regular prism tilings (simply truncated Coxeter orthoscheme tilings) and the corresponding optimal hyperball packings in
$\mathbb{H}^n$ $(n=3,4)$ and we extended the method -- developed in the former paper \cite{Sz06-2} -- to 
$5$-dimensional hyperbolic space (see \cite{Sz13-3}).
In paper \cite{Sz13-4} we studied the $n$-dimensional hyperbolic regular prism honeycombs
and the corresponding coverings by congruent hyperballs and we determined their least dense covering densities.
Furthermore, we formulated conjectures for candidates of the least dense hyperball
covering by congruent hyperballs in $3$- and $5$-dimensional hyperbolic spaces.

In \cite{Sz17-1} we discussed congruent and non-congruent hyperball packings of the truncated regular tetrahedron tilings.
These are derived from the Coxeter simplex tilings $\{p,3,3\}$ $(7\le p \in \mathbb{N})$ and $\{5,3,3,3,3\}$
in $3$- and $5$-dimensional hyperbolic space.
We determined the densest hyperball packing arrangement and its density
with congruent hyperballs in $\mathbb{H}^5$ and determined the smallest density upper bounds of 
non-congruent hyperball packings generated by the above tilings in $\HYN,~ (n=3,5)$.

In \cite{Sz17} we deal with the packings derived by horo- and hyperballs (briefly hyp-hor packings) in $n$-dimensional hyperbolic spaces $\HYN$
($n=2,3$) which form a new class of the classical packing problems.
We constructed in the $2-$ and $3-$dimensional hyperbolic spaces hyp-hor packings that
are generated by complete Coxeter tilings of degree $1$ 
and we determined their densest packing configurations and their densities.
We proved using also numerical approximation methods that in the hyperbolic plane ($n=2$) the density of the above hyp-hor packings arbitrarily approximate
the universal upper bound of the hypercycle or horocycle packing density $\frac{3}{\pi}$ and
in $\HYP$ the optimal configuration belongs to the $\{7,3,6\}$ Coxeter tiling with density $\approx 0.83267$.
Furthermore, we analyzed the hyp-hor packings in
truncated orthosche\-mes $\{p,3,6\}$ $(6< p < 7, ~ p\in \mathbb{R})$ whose
density function is attained its maximum for a parameter which lies in the interval $[6.05,6.06]$
and the densities for parameters lying in this interval are larger that $\approx 0.85397$. 

In \cite{Sz14} we proved that if the truncated tetrahedron is regular, then the density
of the densest packing is $\approx 0.86338$. This is larger than the B\"or\"oczky-Florian density upper bound
but our locally optimal hyperball packing configuration cannot be extended to the entirety of
$\mathbb{H}^3$. However, we described a hyperball packing construction, 
by the regular truncated tetrahedron tiling under the extended Coxeter group $\{3, 3, 7\}$ with maximal density $\approx 0.82251$.

Recently, (to the best of author's knowledge) the candidates for the densest hyperball
(hypersphere) packings in the $3,4$ and $5$-dimensional hyperbolic space $\mathbb{H}^n$ are derived by the regular prism
tilings which have been in papers \cite{Sz06-1}, \cite{Sz06-2} and \cite{Sz13-3}.

{In \cite{Sz17-2} we considered hyperball packings in 
$3$-dimensional hyperbolic space. We developed a decomposition algorithm that for each saturated hyperball packing provides a decomposition of $\HYP$ 
into truncated tetrahedra. Therefore, in order to get a density upper bound for hyperball packings, it is sufficient to determine
the density upper bound of hyperball packings in truncated simplices.}

In \cite{Sz18} we studied hyperball packings related to truncated regular octahedron and cube tilings that are derived from the Coxeter simplex tilings 
$\{p,3,4\}$ $(7\le p \in \mathbb{N})$ and $\{p,4,3\}$ $(5\le p \in \mathbb{N})$
in $3$-dimensional hyperbolic space $\HYP$. We determined the densest hyperball packing arrangement and its density
with congruent and non-congruent hyperballs related to the above tilings. 
Moreover, we prove that the locally densest congruent or non-congruent hyperball configuration belongs to the regular truncated cube with density
$\approx 0.86145$. This is larger than the B\"or\"oczky-Florian density upper bound for balls and horoballs.
We described a non-congruent hyperball packing construction, by the regular cube tiling under the extended Coxeter group $\{4, 3, 7\}$ 
with maximal density $\approx 0.84931$.

{\it In the present paper we study congruent and non-congruent hyperball packings generated by doubly truncated Coxeter orthoscheme tilings in the $3$-dimensional hyperbolic space.
We prove that the densest congruent hyperball packing belongs to the Coxeter orthoscheme tiling of parameter $\{7,3,7\}$ with density $\approx 0.81335$ (see Theorems 4.4-5).
This density is equal -- by our conjecture -- with the upper bound density of the corresponding non-congruent hyperball arrangements (see Theorem 4.8 and Conjecture 4.1).   
}
\section{Orthoschemes, hyperspheres and their volumes}
{\it An orthoscheme $\mathcal{O}$ in $\mathbb{H}^n$ $n \geq 2$ in classical sense}
is a simplex bounded by $n+1$ hyperplanes $H_0,\dots,H_n$
such that $H_i \bot H_j, \  \text{for} \ j \ne i-1,i,i+1.$ Or,
equivalently, the $n+1$ vertices of $\mathcal{O}$ can be
labelled by $A_0,A_1,\dots,A_n$ in such a way that
$\text{span}(A_0,\dots,A_i) \perp \text{span}(A_i,\dots,A_n) \ \ \text{for} \ \ 0<i<n-1.$

Geometrically, complete orthoschemes of degree $m=0,1,2$ can be described as
follows:

\begin{enumerate}
\item
For $m=0$, they coincide with the class of classical orthoschemes introduced by
{{Schl\"afli}}. The initial and final vertices, $A_0$ and $A_n$ of the orthogonal edge-path
$A_iA_{i+1},~ i=0,\dots,n-1$, are called principal vertices of the orthoscheme (see Remark 4.1).
\item
A complete orthoscheme of degree $m=1$ can be constructed from an
orthoscheme with one outer principal vertex, one of $A_0$ or $A_n$, which is simply truncated by
its polar plane (see Fig.~1-2).
\item
A complete orthoscheme of degree $m=2$ can be constructed from an
orthoscheme with two outer principal vertices, $A_0$ and $A_n$, which is doubly truncated by
their polar planes $pol(A_0)$ and $pol(A_n)$ (see Fig.~1-2).
\end{enumerate}
For the {\it complete Coxeter orthoschemes} $\mathcal{O} \subset \mathbb{H}^n$ we adopt the usual
conventions and sometimes even use them in the Coxeter case: If two nodes are related by the weight $\cos{\frac{\pi}{p}}$
then they are joined by a ($p-2$)-fold line for $p=3,~4$ and by a single line marked by $p$ for $p \geq 5$.
In the hyperbolic case if two bounding hyperplanes of $S$ are parallel, then the corresponding nodes
are joined by a line marked $\infty$. If they are divergent then their nodes are joined by a dotted line.

In the following we concentrate only on dimensions $3$ and on hyperbolic
Coxeter-Schl\"afli symbol of the complete orthoscheme tiling $\mathcal{P}$ generated by reflections in the planes of a complete orthoscheme $\mathcal{O}$.
To every scheme there is a corresponding
symmetric $4 \times 4$ matrix $(b^{ij})$ where $b^{ii}=1$ and, for $i \ne j\in \{0,1,2,3\}$,
$b^{ij}$ equals $-\cos{\alpha_{ij}}$ with all angles $\alpha_{ij}$ between the faces $i$,$j$ of $\mathcal{O}$.

For example, $(b^{ij})$ below is the so called Coxeter-Schl\"afli matrix with
parameters $(u;v;w)$, i.e. $\alpha_{01}=\frac{\pi}{u}$, $\alpha_{12}=\frac{\pi}{v}$, $\alpha_{23}=\frac{\pi}{w}$
to be discussed yet for hyperbolicity. Now only $3\le u,v,w$ come into account (see \cite{IH90}). Then it holds
\[
(b^{ij})=\langle \mbox{\boldmath$b^i$},\mbox{\boldmath$b^j$} \rangle:=\begin{pmatrix}
1& -\cos{\frac{\pi}{u}}& 0 & 0 \\
-\cos{\frac{\pi}{u}} & 1 & -\cos{\frac{\pi}{v}}& 0 \\
0 & -\cos{\frac{\pi}{v}} & 1 & -\cos{\frac{\pi}{w}} \\
0 & 0 & -\cos{\frac{\pi}{w}} & 1
\end{pmatrix}. \tag{2.1}
\]

This $3$-dimensional complete (truncated or frustum) orthoscheme $\mathcal{O}=W_{uvw}$ and its reflection group $\bG_{uvw}$ will be described in
Fig.~2, and by the symmetric Coxeter-Schl\"afli matrix $(b^{ij})$ in formula (2.1), furthermore by its inverse matrix $(a_{ij})$ in formula
(2.2).
\[
\begin{gathered}
(a_{ij})=(b^{ij})^{-1}=\langle \ba_i, \ba_j \rangle:=\\
=\frac{1}{B} \begin{pmatrix}
\sin^2{\frac{\pi}{w}}-\cos^2{\frac{\pi}{v}}& \cos{\frac{\pi}{u}}\sin^2{\frac{\pi}{w}}& \cos{\frac{\pi}{u}}\cos{\frac{\pi}{v}} & \cos{\frac{\pi}{u}}\cos{\frac{\pi}{v}}\cos{\frac{\pi}{w}} \\
\cos{\frac{\pi}{u}}\sin^2{\frac{\pi}{w}} & \sin^2{\frac{\pi}{w}} & \cos{\frac{\pi}{v}}& \cos{\frac{\pi}{w}}\cos{\frac{\pi}{v}} \\
\cos{\frac{\pi}{u}}\cos{\frac{\pi}{v}} & \cos{\frac{\pi}{v}} & \sin^2{\frac{\pi}{u}}  & \cos{\frac{\pi}{w}}\sin^2{\frac{\pi}{u}}  \\
\cos{\frac{\pi}{u}}\cos{\frac{\pi}{v}}\cos{\frac{\pi}{w}}  & \cos{\frac{\pi}{w}}\cos{\frac{\pi}{v}} & \cos{\frac{\pi}{w}}\sin^2{\frac{\pi}{u}}  & \sin^2{\frac{\pi}{u}}-\cos^2{\frac{\pi}{v}}
\end{pmatrix}, \tag{2.2}
\end{gathered}
\]
where
$$
B=\det(b^{ij})=\sin^2{\frac{\pi}{u}}\sin^2{\frac{\pi}{w}}-\cos^2{\frac{\pi}{v}} <0, \ \ \text{i.e.} \ \sin{\frac{\pi}{u}}\sin{\frac{\pi}{w}}-\cos{\frac{\pi}{v}}<0.
$$

In the following we use the above orthoscheme whose volume is derived by the next
Theorem of {{R.~Kellerhals}} (\cite{K89}, by the ideas of N.~I.~Lobachevsky):
\begin{theorem}{\rm{(R.~Kellerhals)}} The volume of a three-dimensional hyperbolic
complete orthoscheme $\mathcal{O}=W_{uvw} \subset \mathbb{H}^3$
is expressed with the essential
angles $\alpha_{01}=\frac{\pi}{u}$, $\alpha_{12}=\frac{\pi}{v}$, $\alpha_{23}=\frac{\pi}{w}$, $(0 \le \alpha_{ij}
\le \frac{\pi}{2})$
(Fig.~1.) in the following form:

\begin{align}
&\mathrm{Vol}(\mathcal{O})=\frac{1}{4} \{ \mathcal{L}(\alpha_{01}+\theta)-
\mathcal{L}(\alpha_{01}-\theta)+\mathcal{L}(\frac{\pi}{2}+\alpha_{12}-\theta)+ \notag \\
&+\mathcal{L}(\frac{\pi}{2}-\alpha_{12}-\theta)+\mathcal{L}(\alpha_{23}+\theta)-
\mathcal{L}(\alpha_{23}-\theta)+2\mathcal{L}(\frac{\pi}{2}-\theta) \}, \notag
\end{align}
where $\theta \in [0,\frac{\pi}{2})$ is defined by:
$$
\tan(\theta)=\frac{\sqrt{ \cos^2{\alpha_{12}}-\sin^2{\alpha_{01}} \sin^2{\alpha_{23}
}}} {\cos{\alpha_{01}}\cos{\alpha_{23}}},
$$
and where $\mathcal{L}(x):=-\int\limits_0^x \log \vert {2\sin{t}} \vert dt$ \ denotes the
Lobachevsky function (in J. Milnor's interpretation).
\end{theorem}
The hypersphere (or equidistant surface) is a quadratic surface at a constant distance
from a plane (base plane) in both halfspaces. The infinite body of the hypersphere, containing the base plane, is called {\it hyperball}.

The {\it half hyperball } (i.e., the part of the hyperball lying on one side of its base plane) with distance $h$ to a base plane $\beta$ is denoted by $\mathcal{H}^h_+$.
The volume of the intersection of $\mathcal{H}^h_+(\mathcal{A})$ and the right prism with
base a $2$-polygon $\mathcal{A} \subset \beta$, can be determined by the classical formula 
(2.1) of J.~Bolyai \cite{B31}.
\begin{equation}
\mathrm{Vol}(\mathcal{H}^h_+(\mathcal{A}))=\frac{1}{4}\mathrm{Area}(\mathcal{A})\left[k \sinh \frac{2h}{k}+
2 h \right], \tag{2.3}
\end{equation}
The constant $k =\sqrt{\frac{-1}{K}}$ is the natural length unit in
$\mathbb{H}^3$, where $K$ denotes the constant negative sectional curvature. In the following we may assume that $k=1$.
\section{Essential points in a doubly truncated orthoscheme}
Let $A_0(\ba_0)$, $A_1(\ba_1)$, $A_2(\ba_2)$, $A_3(\ba_3)$ be the vertices of the above complete orthoscheme $W_{uvw}$ (see Fig.~1,2).
In the considered cases the principal vertices $A_0$ and $A_3$ are outer points $(a_{ii}>0)$, ($i \in \{0,3\}$). 

We distinguish the following main configurations:
\begin{enumerate}
\item[{\bf 1.}] {\it $A_3$ is outer point $\frac{\pi}{u}+\frac{\pi}{v} < \frac{\pi}{2}$, then $a_3(\mbox{\boldmath$a$}_{3})=JEQ$ is its polar plane and 
$A_0$ is also outer $\frac{\pi}{v}+\frac{\pi}{w} < \frac{\pi}{2}$, then $a_0(\mbox{\boldmath$a$}_{0})=CLH$ is its polar plane.}
\item[{\bf 1.i}] $u=w$, $F_{03}F_{12}$ is half turn axis, $\mathbf{h}$ is the half turn changing $0 \leftrightarrow 3$, $1 \leftrightarrow 2$.
Here a "half orthoscheme" $JQEB_{13}F_{12}B_{02}F_{03}A_2$ will be the fundamental domain of $\bG_{u=w,v}$.
\end{enumerate}
We use for calculations the following important lemmas, (see \cite{MoSzi18} and Fig.~1,~2):
\begin{lemma}
Let $A_0$ be an outer principal vertex of the orthoscheme $W_{uvw}$ and let $a_0(\mbox{\boldmath$a$}_{0})=CLH$ be its polar plane where
$C=a_0 \cap A_0A_1$, $L=a_0 \cap A_0A_2$, $H=a_0 \cap A_0A_3$ whose vectors are the following:
\begin{equation}
\begin{gathered}
C(\bc)=a_0 \cap A_0A_1; ~ \bc=\ba_1-\frac{a_{01}}{a_{00}} \ba_0, \ \mathrm{with} \\
\langle \bc,\bc \rangle =\frac{(a_{11}a_{00}-a_{01}^2)}{a_{00}}=\langle \bc,\ba_1 \rangle
=\frac{\sin^2\frac{\pi}{w}}{\sin^2\frac{\pi}{w}-\cos^2\frac{\pi}{v}}=\frac{a_{11}}{a_{00}} \\
L(\mathbf{l})=a_0 \cap A_0A_2; ~ \mathbf{l}=\ba_2-\frac{a_{02}}{a_{00}} \ba_0, \ \mathrm{with} \\
\langle \bl,\bl \rangle =\frac{(a_{22}a_{00}-a_{02}^2)}{a_{00}}=\langle \bl,\ba_2 \rangle
=\frac{1}{\sin^2\frac{\pi}{w}-\cos^2\frac{\pi}{v}}=\frac{1}{Ba_{00}} \\
H(\bh)=a_0 \cap A_0A_3; ~ \bh=\ba_3-\frac{a_{03}}{a_{00}} \ba_0, \ \mathrm{with} \\
\langle \bh,\bh \rangle =\frac{(a_{33}a_{00}-a_{03}^2)}{a_{00}}=\langle \bh,\ba_3 \rangle
=\frac{\sin^2\frac{\pi}{v}}{\sin^2\frac{\pi}{w}-\cos^2\frac{\pi}{v}}=\frac{\sin^2\frac{\pi}{v}}{Ba_{00}}
\end{gathered} \tag{3.1}
\end{equation}
\end{lemma}
\begin{figure}[htb]
\centering
\includegraphics[width=80mm]{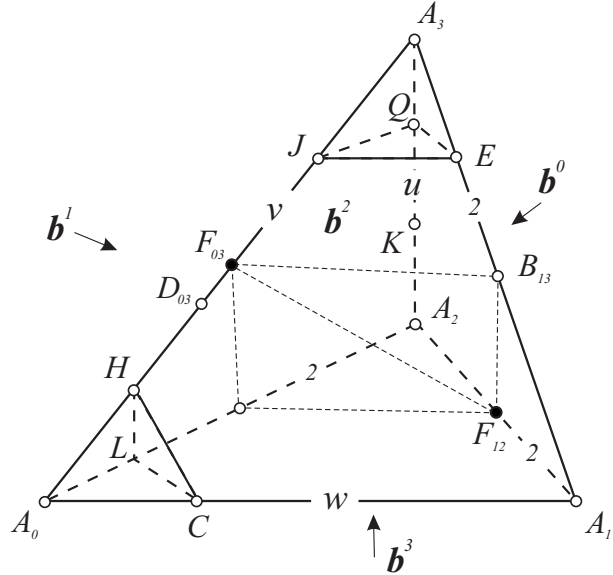} 
\caption{Double truncated complete orthoscheme with essential points }
\label{}
\end{figure}
\begin{lemma}
Let $A_3$ be an outer principal vertex of the orthoscheme $W_{uvw}$ and let $a_3(\mbox{\boldmath$a$}_{3})=JEQ$ be its polar plane where
$J=a_3 \cap A_3A_0$, $E=a_3 \cap A_3A_1$, $Q=a_3 \cap A_3A_2$ whose vectors are the following:
\begin{equation}
\begin{gathered}
J(\bj)=a_3 \cap A_3A_0; ~ \bj=\ba_0-\frac{a_{03}}{a_{33}} \ba_3, \ \mathrm{with} \\
\langle \bj,\bj \rangle =\frac{(a_{00}a_{33}-a_{03}^2)}{a_{33}}=\langle \bj,\ba_0 \rangle
=\frac{\sin^2\frac{\pi}{v}}{\sin^2\frac{\pi}{u}-\cos^2\frac{\pi}{v}}=\frac{\sin^2\frac{\pi}{v}}{Ba_{33}} \\
E(\mathbf{e})=a_3 \cap A_3A_1; ~ \mathbf{e}=\ba_1-\frac{a_{13}}{a_{33}} \ba_3, \ \mathrm{with} \\
\langle \be,\be \rangle =\frac{(a_{11}a_{33}-a_{13}^2)}{a_{33}}=\langle \be,\ba_1 \rangle
=\frac{1}{\sin^2\frac{\pi}{u}-\cos^2\frac{\pi}{v}}=\frac{1}{Ba_{33}} \\
Q(\bh)=a_3 \cap A_3A_2; ~ \bq=\ba_2-\frac{a_{23}}{a_{33}} \ba_3, \ \mathrm{with} \\
\langle \bq,\bq \rangle =\frac{(a_{22}a_{33}-a_{23}^2)}{a_{33}}=\langle \bq,\ba_2 \rangle
=\frac{\sin^2\frac{\pi}{u}}{\sin^2\frac{\pi}{u}-\cos^2\frac{\pi}{v}}=\frac{a_{22}}{a_{33}}.
\end{gathered} \tag{3.2}
\end{equation}
\end{lemma}
\begin{figure}[htb]
\centering
\includegraphics[width=90mm]{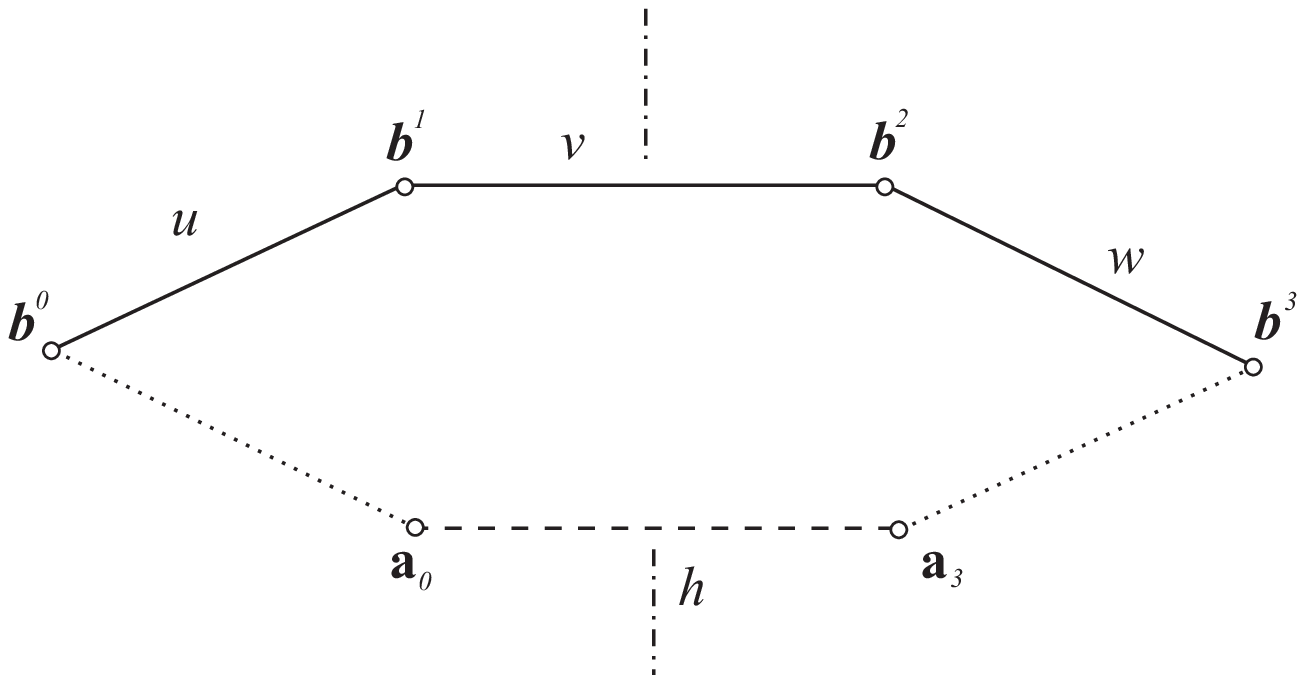}
\caption{}
\label{Fig2}
\end{figure}
Especially if $u=w$, the midpoints $F_{03}$ of $JH$ and $F_{12}$ of $A_1A_2$ can play important roles, since
$F_{03}F_{12}$ will be the axis of half turn
$$
\bh:~ 0 \leftrightarrow 3, 1 \leftrightarrow 2, ~\mathrm{i.e.} ~ A_0 \leftrightarrow A_3, ~
\mbox{\boldmath$b$}^{0} \leftrightarrow \mbox{\boldmath$b$}^{3},  ~ A_1 \leftrightarrow A_2, ~
\mbox{\boldmath$b$}^{1} \leftrightarrow \mbox{\boldmath$b$}^{2}.
$$
(Here $a_{00}=a_{33}$ and $a_{11}=a_{22}$ hold, of course.)
\begin{lemma}
The midpoints $F_{03}(\mathbf{f}_{03})$ of $JH$ and  $F_{12}(\mathbf{f}_{12})$ of $A_1 A_2$ (see Fig.~1) can be
determined by the following vectors:
\begin{equation}
\begin{gathered}
\mathbf{f}_{03}=\ba_0+\ba_3, ~ \langle \mathbf{f}_{03},\mathbf{f}_{03} \rangle= 2(a_{00}+a_{03})<0, \\
\mathbf{f}_{12}=\ba_1+\ba_2, ~ \langle \mathbf{f}_{12},\mathbf{f}_{12} \rangle= 2(a_{11}+a_{12})<0.
\end{gathered} \notag
\end{equation}
\end{lemma}
\section{On hyperball packings in a doubly truncated orthoscheme}
Similarly to the former cases (see \cite{Sz06-1}, \cite{Sz06-2}, \cite{Sz13-4}, \cite{Sz14}, \cite{Sz17-1}, \cite{Sz17-2}) it is interesting to study and 
to construct locally optimal {\it congruent and non-congruent} hyperball packings relating to suitable truncated polyhedron tilings in $3$- and higher dimensions as well.
This study fits into our program to look for the upper bound density of the congruent and non-congruent hyperball packings in $\HYN$.

{\it In this paper we consider the 3-dimensional regular doubly truncated orthoscheme tilings and study their corresponding locally optimal congruent and non-congruent hyperball
packings.}

\subsection{Congruent hyperball packings}
We consider a doubly truncated orthoscheme tiling $\mathcal{T}(\mathcal{O}(u,v,w))$ with Schl\"afli symbol $\{u,v,w\}$, 
($\frac{1}{u}+\frac{1}{v} < \frac{1}{2},~\frac{1}{v}+\frac{1}{w} < \frac{1}{2}$, $3 \le u,v,w \in \mathbb{N}$) whose fundamental 
domain is doubly truncated orthoschem (e.g. $CHLA_1A_2EJQ$ in Fig.~1). 

Let a truncated orthoscheme $\mathcal{O}(u,v,w)$ $\subset \HYP$ be a tile from 
the above tiling. 
This truncated orthoscheme can be derived also by truncation from 
a orthoschem $A_0A_1A_2A_3$ with outer essential vertices $A_0$ and $A_3$.
The truncating planes $a_0(\mbox{\boldmath$a$}_{0})=CLH$ and $a_3(\mbox{\boldmath$a$}_{3})=JEQ$ are the polar planes of outer vertices $A_0$ and $A_3$, 
that {\it can be the ultraparallel base planes 
of hyperballs $\mathcal{H}^{s}_i$} with height $s$ $(i=0,3)$. 
The distance between the two base planes is $2h^{03}(u,v,w)=d(a_0(\mbox{\boldmath$a$}_{0}),a_3(\mbox{\boldmath$a$}_{3}))=d(H,J)$ ($d$ is the hyperbolic distance function).
\begin{figure}[htb]
\centering
\includegraphics[width=120mm]{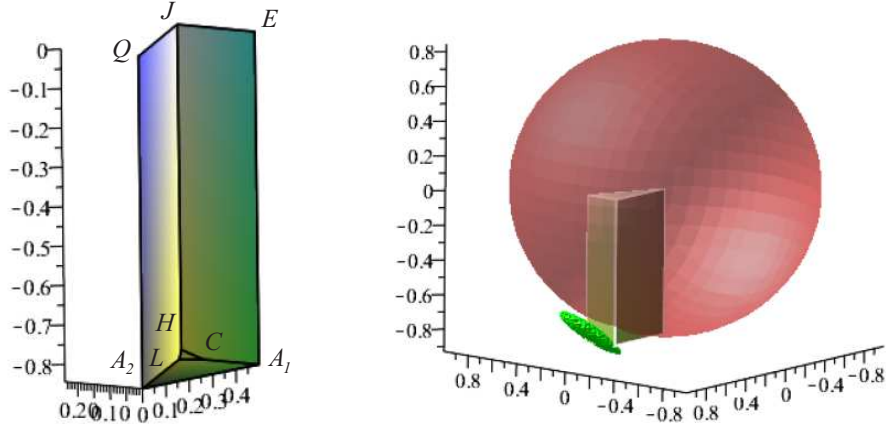} 
\caption{The densest congruent hyperball packing arrangement related to parameters $\{7,3,7\}$ with density $\approx 0.81335$ }
\label{}
\end{figure}
In this subsection we consider congruent hyperball packings therefore we have to distinguish $3$ different cases.
\begin{enumerate}
\item Both polar planes are assigned hyperspheres that are congruent with each other therefore the height of a hyperball is at most $h^{03}(u,v,w)$ (see Fig.~1). It is clear, 
that the heights $h^{0}=h^3$ of optimal hyperballs $\mathcal{H}_i^{h^{i}}$ $(i=0,3)$ is
\begin{equation}
\begin{gathered}
h=h^{0}(u,v,w)=h^3(u,v,w)=\\ 
=\min\{h^{03}(u,v,w)=d(H,J)/2,d(Q,A_2),d(C,A_1)\}, \tag{4.1}
\end{gathered}
\end{equation}
where $u,v,w$ are suitable given integer parameters. 
In this case the volume sum of the hyperball pieces 
lying in the orthoscheme is 
$$Vol(\mathcal{H}^{h}(\mathcal{A}_0)\cap \mathcal{O}(u,v,w))+
Vol(\mathcal{H}^{h}(\mathcal{A}_3)\cap \mathcal{O}(u,v,w))$$ (see (2.3))
where $\mathcal{A}_0$ is the area of the triangle $CLH$ and 
$\mathcal{A}_3$ is the area of the triangle $JEQ$. 

\item In these cases we consider only one hyperball type:
\begin{enumerate}
\item with base plane $a_0(\mbox{\boldmath$a$}_{0})=CLH$.  
The height of the optimal hyperball $\mathcal{H}_0^{h^0}$ is
\begin{equation}
h^0(u,v,w)=\min\{2h^{03}(u,v,w)=d(J,H),d(C,A_1)\}, \tag{4.2}
\end{equation}
where $u,v,w$ are suitable given integer parameters.
\item with base plane $a_3(\mbox{\boldmath$a$}_{3})=JEQ$.
The height of the optimal hyperball $\mathcal{H}_3^{h^3}$ is
\begin{equation}
h^3(u,v,w)=\min\{2h^{03}(u,v,w)=d(J,H),d(Q,A_2)\}, \tag{4.3}
\end{equation}
where $u,v,w$ are suitable given integer parameters.
\end{enumerate}
\end{enumerate}
The volume of a doubly truncated orthoscheme $\mathcal{O}(u,v,w)$ is denoted by $Vol(\mathcal{O}(u,v,w))$.
\begin{definition} The locally density functions $\delta^i(\mathcal{O}(u,v,w)$ of the congruent hyperball packings related to 
$\mathcal{O}(u,v,w)$ and the above cases ($i\in \{1,2\}$) are 
defined by next formulas:
\begin{enumerate}
\item 
\begin{equation}
\delta^1(\mathcal{O}(u,v,w)):=\frac{Vol(\mathcal{H}^{s}(\mathcal{A}_0)\cap \mathcal{O}(u,v,w))+
Vol(\mathcal{H}_3^{s}(\mathcal{A}_3)\cap \mathcal{O}(u,v,w))}{Vol({\mathcal{O}(u,v,w)})}, \notag
\end{equation}
where $0<s\le h$, $\mathcal{A}_0$ is the area of the triangle $CLH$ and 
$\mathcal{A}_3$ is the area of the triangle $JEQ$ (see (4.1)).
\item
\begin{equation}
\delta^2_j(\mathcal{O}(u,v,w)):=\frac{Vol(\mathcal{H}^{s}(\mathcal{A}_j) \cap \mathcal{O}(u,v,w))}{Vol({\mathcal{O}(u,v,w)})}, \notag
\end{equation}
where $0<s\le h^{j} $, $j\in\{0,3\}$ (see (4.2), (4.3)) and
\begin{equation}
\delta^2(\mathcal{O}(u,v,w)):=\max_{j=0,3}\{\delta^2_j(\mathcal{O}(u,v,w))\}. \notag
\end{equation}
\end{enumerate}
\end{definition}
The distance $s$ of two proper points
$X(\mathbf{x})$ and $Y(\mathbf{y})$ is calculated by the formula
\begin{equation}
\cosh{{s}}=\frac{-\langle ~ \mathbf{x},~\mathbf{y} \rangle }{\sqrt{\langle ~ \mathbf{x},~\mathbf{x} \rangle
\langle ~ \mathbf{y},~\mathbf{y} \rangle }}. \tag{4.4}
\end{equation}
If the parameters $u,v,w$ are given then the lengths of the line segments $A_1C$, $A_2Q$ and $JH$ 
can be determined by the machinery of the projective geometry using the Lemmas 3.1-3 and formula (4.4):
\begin{lemma}
\begin{equation}
\begin{gathered}
d(A_1,C)=\mathrm{arcosh}{\frac{1}{\sqrt{a_{00}}}},~ 
d(A_2,Q)=\mathrm{arcosh}{\frac{1}{\sqrt{a_{33}}}},~ d(J,H)=\mathrm{arcosh}{\frac{-a_{03}}{\sqrt{a_{00}a_{33}}}}.
\end{gathered} \notag 
\end{equation} 
where $(a_{ij})$ ($i,j=0,1,2,3)$) is the inverse of the corresponding Coxeter-Schl\"afli matrix (see (2.2).
\end{lemma}

In our cases the essential dihedral angles of orthoschemes $\mathcal{O}(u,v,w)$ 
are the following:
$\alpha_{01}=\frac{\pi}{u}, \ \ \alpha_{12}=\frac{\pi}{v}, \ \
\alpha_{23}=\frac{\pi}{w}$ (see Fig.~1), therefore, the volume 
${Vol}(\mathcal{O}(u,v,w))$ of the orthoscheme $\mathcal{O}(u,v,w)$ can be determined by Theorem 2.1.
Moreover, the maximal height $h$, $h^0$ or $h^3$ of congruent optimal hyperballs 
and the the corresponding volumes of the hyperball pieces  can be computed for any suitable 
fixed integer parameters $u,v,w$. Therefore, the density 
$\delta^1(\mathcal{O}(u,v,w))$ or $\delta^2_j(\mathcal{O}(u,v,w))$ $(j\in\{1,2\})$ (see Definition 4.1)
depends only on the suitable integer parameters $u,v,w$ of the doubly truncated orthoscheme $\mathcal{O}(u,v,w)$.
Moreover, the volumes of the hyperball pieces can be computed by the formula (2.3).
\subsubsection{Numerical data of the optimal congruent hyperball arrangement}
First, we illustrate our computation method for given important parameters and 
then we summarize the numerical data of optimal congruent hyperball arrangements for 
several parameters in Tables 1 and 2.
\medbreak
{\bf Results for parameters $u=7$, $v=3$, $w=7$}:

\begin{equation}
\begin{gathered}
d(A_1,C)=\mathrm{arcosh}{\frac{1}{\sqrt{a_{00}}}}= d(A_2,Q)=\mathrm{arcosh}{\frac{1}{\sqrt{a_{33}}}} \approx 1.23469,\\
d(J,H)/2=\frac{1}{2}\mathrm{arcosh}{\frac{-a_{03}}{\sqrt{a_{00}a_{33}}}}\approx 1.28517
\end{gathered} \notag 
\end{equation}
Therefore, the optimal heights in all cases are equal: $h=h^0=h^3=d(A_1,C)=d(A_2,Q) \approx 1.23469$ (see (4.1), (4.2), (4.3)). 
\begin{equation}
\begin{gathered}
Vol({\mathcal{O}(7,3,7)}) \approx 0.38325,~ ~ 
Vol(\mathcal{H}^{h}(\mathcal{A}_0))=Vol(\mathcal{H}^{h}(\mathcal{A}_3))=\\ 
=Vol(\mathcal{H}^{h^0}(\mathcal{A}_0))=Vol(\mathcal{H}^{h^3}(\mathcal{A}_3)) \approx 0.15586,
\end{gathered} \notag
\end{equation} 
\begin{enumerate}
\item {\it Two congruent hyperball:}
\begin{equation} 
\begin{gathered}
\delta^1(\mathcal{O}(7,3,7))=\frac{2 \cdot Vol(\mathcal{H}^{h}(\mathcal{A}_0) \cap \mathcal{O}(7,3,7))}{Vol({\mathcal{O}(7,3,7)})}\approx 0.81335,
\end{gathered} \notag
\end{equation}
\item {\it One hyperball:}
\begin{equation} 
\begin{gathered}
\delta^2_j(\mathcal{O}(7,3,7))=\frac{Vol(\mathcal{H}^{h}(\mathcal{A}_j) \cap \mathcal{O}(7,3,7))}{Vol({\mathcal{O}(7,3,7)})} \approx 0.40668, (j\in \{0,3\})
\end{gathered} \notag
\end{equation}
\end{enumerate}
\begin{rmrk}
If $u=w$ and $h \ge h^0=h^3$ then $\delta^2(\mathcal{O}(u,v,w))=
\frac{1}{2}\delta^1(\mathcal{O}(u,v,w))$ (see above example with parameters $u=7,v=3,w=7$). 
\end{rmrk}
In the following Table we summarize the data of the hyperball packings for some 
parameters $u,v,w \in \mathbb{N}$, where $\mathcal{A}_i$ ($i\in \{0,3\}$) 
is the area of the trigonal face (triangle $CHL$ or $EJQ$) of the truncated 
tetrahedron related to the vertex $A_i$, cf. Fig.~1 (see (4.1), (4.2), (4.3) and Definition 4.1).
We note here, that the role of the parameters $u$ and $w$ is symmetrical 
therefore we can assume, that $u \le w$.
\medbreak
{\scriptsize
\centerline{\vbox{
\halign{\strut\vrule~\hfil $#$ \hfil~\vrule
&\quad \hfil $#$ \hfil~\vrule
&\quad \hfil $#$ \hfil\quad\vrule
&\quad \hfil $#$ \hfil\quad\vrule
&\quad \hfil $#$ \hfil\quad\vrule
\cr
\noalign{\hrule}
\multispan5{\strut\vrule\hfill\bf Table 1, two congruent hyperballs  \hfill\vrule}%
\cr
\noalign{\hrule}
\noalign{\vskip2pt}
\noalign{\hrule}
\{u,v,w\} & h  & Vol({\mathcal{O}(u,v,w)}) & \sum_{i=0,3} {Vol}(\mathcal{H}^h(\mathcal{A}_i))\\
& \delta^1({\mathcal{O}(u,v,w)})\cr
\noalign{\hrule}
\{7,3,7\} & 1.23469 & 0.38325 & 0.31172 & 0.81335 \cr
\noalign{\hrule}
\{7,3,8\} & 0.93100 & 0.41326 & 0.25726 & 0.62251  \cr
\noalign{\hrule}
\{7,3,9\} & 0.76734 & 0.43171 & 0.23355 & 0.54099 \cr
\noalign{\hrule}
\vdots & \vdots  & \vdots  & \vdots  & \vdots  \cr
\noalign{\hrule}
\{7,3,50\} & 0.11380 & 0.49016 & 0.06121 & 0.12488 \cr
\noalign{\hrule}
\vdots & \vdots  & \vdots  & \vdots  & \vdots  \cr
\noalign{\hrule}
\{8,3,8\} & 0.94946 & 0.44383 & 0.33794 & 0.76143 \cr
\noalign{\hrule}
\{8,3,9\} & 0.78366 & 0.46266 & 0.29474 & 0.63704 \cr
\noalign{\hrule}
\{8,3,10\} & 0.67409& 0.47536 & 0.26747 & 0.56266 \cr
\noalign{\hrule}
\vdots & \vdots  & \vdots  & \vdots  & \vdots  \cr
\noalign{\hrule}
\{8,3,50\} & 0.11668 & 0.52248 & 0.06935 & 0.13274 \cr
\noalign{\hrule}
\vdots & \vdots  & \vdots  & \vdots  & \vdots  \cr
\noalign{\hrule}
\{5,4,5\} & 0.88055 & 0.46190 & 0.36007 & 0.77955 \cr
\noalign{\hrule}
\{5,4,6\} & 0.73969 & 0.50747 & 0.37287 & 0.73476 \cr
\noalign{\hrule}
\{5,4,7\} & 0.59326 & 0.53230 & 0.32974 & 0.61947 \cr
\noalign{\hrule}
\vdots & \vdots  & \vdots  & \vdots  & \vdots  \cr
\noalign{\hrule}
\{5,4,50\} & 0.07206 & 0.59291 & 0.06350 & 0.10710 \cr
\noalign{\hrule}
\vdots & \vdots  & \vdots  & \vdots  & \vdots  \cr
\noalign{\hrule}
\{4,5,4\} & 0.80846 & 0.43062 & 0.31702 & 0.73620 \cr
\noalign{\hrule}
\{4,5,5\} & 0.69129 & 0.49789 & 0.38284 & 0.76893 \cr
\noalign{\hrule}
\{4,5,6\} & 0.53064 & 0.52971 & 0.33597 & 0.63426 \cr
\noalign{\hrule}
\vdots & \vdots  & \vdots  & \vdots  & \vdots  \cr
\noalign{\hrule}
\{4,5,50\} & 0.05502 & 0.59318 & 0.05710 & 0.096256 \cr
\noalign{\hrule}}}}}
\medbreak
The volume $Vol({\mathcal{O}(u,v,w)})$ can be calculated by Theorem 2.1. The maximal volume sum $\sum_{i=0,3} {Vol}(\mathcal{H}^h(\mathcal{A}_i))$ of the hyperball 
pieces lying in ${\mathcal{O}(u,v,w)}$ can be computed by the formulas (2.3), (4.1), (4.2), (4.3) and by the above described computation method 
for each given possible parameters $u,v,w$. Therefore, the maximal density of the congruent hyperball packing related to every doubly truncated orthoscheme tiling with ''two hyperballs" 
-- $\delta^1({\mathcal{O}(u,v,w)})$ (see Definition 4.1) -- can be computed for each possible parameters.

Finally, we obtain after careful analysis of the function the following
\begin{theorem}
The density function $\delta^1({\mathcal{O}(u,v,w)})$, ($\frac{1}{u}+\frac{1}{v} < \frac{1}{2},~\frac{1}{v}+\frac{1}{w} < \frac{1}{2}$. $3 \le u,v,w \in \mathbb{N}$)
attains its maximum at parameters $\{u,v,w\}=\{7,3,7\}$ with density $\delta^1({\mathcal{O}(u,v,w)}) \approx 0.81335$ (see Table 1).
\end{theorem}
\medbreak
{\scriptsize
\centerline{\vbox{
\halign{\strut\vrule~\hfil $#$ \hfil~\vrule
&\quad \hfil $#$ \hfil~\vrule
&\quad \hfil $#$ \hfil\quad\vrule
&\quad \hfil $#$ \hfil\quad\vrule
&\quad \hfil $#$ \hfil\quad\vrule
&\quad \hfil $#$ \hfil\quad\vrule
\cr
\noalign{\hrule}
\multispan6{\strut\vrule\hfill\bf Table 2, one hyperball \hfill\vrule}%
\cr
\noalign{\hrule}
\noalign{\vskip2pt}
\noalign{\hrule}
\{u,v,w\} & h^0 & h^3 & Vol({\mathcal{O}(u,v,w)}) & \begin{gathered} \max_{i=0,3} \{Vol(\mathcal{H}^{h^i}(\mathcal{A}_i)\}
\end{gathered} & \delta^2({\mathcal{O}(u,v,w)})\cr
\noalign{\hrule}
\{7,3,8\} & 0.93100 & 1.25596 & 0.41326 & 0.16371 & 0.39614  \cr
\noalign{\hrule}
\{7,3,9\} & 0.76734 & 1.27042 & 0.43171 & 0.16543 & 0.38320 \cr
\noalign{\hrule}
\vdots & \vdots   & \vdots & \vdots  & \vdots  & \vdots  \cr
\noalign{\hrule}
\{7,3,50\} & 0.11380 & 1.32226 & 0.49016 & 0.18040 & 0.36805 \cr
\noalign{\hrule}
\vdots & \vdots  & \vdots & \vdots  & \vdots  & \vdots  \cr
\noalign{\hrule}
\{8,3,9\} & 0.78366 & 0.96206 & 0.46266 & 0.17265 & 0.37316 \cr
\noalign{\hrule}
\{8,3,10\} & 0.67409 & 0.97104 & 0.47536 & 0.17531 & 0.36879 \cr
\noalign{\hrule}
\vdots & \vdots & \vdots & \vdots  & \vdots  & \vdots  \cr
\noalign{\hrule}
\{8,3,50\} & 0.11668 & 1.00753 & 0.52248 & 0.18650 & 0.35695 \cr
\noalign{\hrule}
\vdots & \vdots  & \vdots & \vdots  & \vdots  & \vdots  \cr
\noalign{\hrule}
\{5,4,5\} & 1.02221 & 1.02221 & 0.46190 & 0.22942 & 0.49668 \cr
\noalign{\hrule}
\{5,4,6\} & 0.73969 & 1.07541 & 0.50747 & 0.25088 & 0.49437 \cr
\noalign{\hrule}
\{5,4,7\} & 0.59326 & 1.10694 & 0.53230 & 0.26448 & 0.49686 \cr
\noalign{\hrule}
\vdots & \vdots  & \vdots &\vdots  & \vdots  & \vdots  \cr
\noalign{\hrule}
\{5,4,50\} & 0.07206 & 1.19054 & 0.59291 & 0.30407 & 0.51284 \cr
\noalign{\hrule}
\vdots & \vdots  & \vdots & \vdots  & \vdots  & \vdots  \cr
\noalign{\hrule}
\{4,5,4\} & 1.06128 & 1.06128 & 0.43062 & 0.24500 & 0.56895 \cr
\noalign{\hrule}
\{4,5,5\} & 0.69129 & 1.16974 & 0.49789 & 0.29371 & 0.58990 \cr
\noalign{\hrule}
\{4,5,6\} & 0.53064 & 1.22646 & 0.52971 & 0.32284 & 0.60946 \cr
\noalign{\hrule}
\vdots & \vdots  & \vdots & \vdots  & \vdots  & \vdots  \cr
\noalign{\hrule}
\{4,5,50\} & 0.05502 & 1.19344 & 0.59318 & 0.30555 & 0.51510 \cr
\noalign{\hrule}
\vdots & \vdots  & \vdots & \vdots  & \vdots  & \vdots  \cr
\noalign{\hrule}
\{4,6,4\} & 0.88137 & 0.88137 & 0.50192 & 0.30049 & 0.59868 \cr
\noalign{\hrule}
\{4,6,5\} & 0.61415 & 0.97970 & 0.55992 & 0.35582 & 0.63548 \cr
\noalign{\hrule}
\{4,6,6\} & 0.48121 & 1.01251 & 0.58850 & 0.32284 & 0.58711 \cr
\noalign{\hrule}
\vdots & \vdots  & \vdots & \vdots  & \vdots  & \vdots  \cr
\noalign{\hrule}
\{4,6,50\} & 0.05138 & 0.88231 & 0.64697 & 0.30100 & 0.46522 \cr
\noalign{\hrule}}}}}
\medbreak
Similarly to the above case the volume of the doubly truncated orthoscheme can be computed by theorem 2.1 but here we apply only one hyperball type. 
Now, the volume -- $\max_{i=0,3} \{Vol(\mathcal{H}^{h^i}(\mathcal{A}_i)\}$ -- of the optimal hyperball 
piece lying in ${\mathcal{O}(u,v,w)}$ can be computed by the formula (2.3), (4.1), (4.2), (4.3) and by the above described computation method 
for each given possible parameters $u,v,w$. Therefore, the maximal density of the congruent hyperball packing related to the doubly truncated orthoscheme tilings with ''one hyperballs" 
-- $\delta^2({\mathcal{O}(u,v,w)})$ (see Definition 4.1) -- can be computed for each possible parameters.

Finally, we obtain after careful analysis of the function the following (see Fig.~3)
\begin{theorem}
The density function $\delta^2({\mathcal{O}(u,v,w)})$, ($\frac{1}{u}+\frac{1}{v} < \frac{1}{2},~\frac{1}{v}+\frac{1}{w} < \frac{1}{2}$. $3 \le u,v,w \in \mathbb{N}$)
attains its maximum at parameters $\{u,v,w\}=\{4,6,5\}$ with density $\delta^2({\mathcal{O}(u,v,w)}) \approx 0.63548$ (see Table 2).
\end{theorem}
\subsection{Non-congruent hyperball packings}
In this subsection we consider non-congruent hyperball packings related to the doubly truncated 
Coxeter orthoschemes which can be derived, similarly to the above section by truncation from 
a orthoschem $A_0A_1A_2A_3$ with outer essential vertices $A_0$ and $A_3$.
The truncating planes are the polar planes of outer vertices $A_0$ and $A_3$, 
that {\it can be the ultraparallel base planes 
of hyperballs $\mathcal{H}^{s}_i$ $(i=0,3)$} with height $s$. 
The distance between the two base planes is $2h^{03}(u,v,w)=d(H,J)$.
\begin{lemma}
If $u \le w$ then $d(C,A_1) \le d(Q,A_2)$ where $u,v,w$ are suitable given integer parameters.
\end{lemma}
{\bf{Proof:}}~ We get from the Lemma 4.2 that $d(A_1,C)=\mathrm{arcosh}{\frac{1}{\sqrt{a_{00}}}}$, $d(A_2,Q)=\mathrm{arcosh}{\frac{1}{\sqrt{a_{33}}}}$  
where $(a_{ij})$ ($i,j=0,1,2,3)$) is the inverse of the corresponding Coxeter-Schl\"afli matrix (see (2.2)). From the matrix $(a_{ij})$ follows that
\begin{equation}
\begin{gathered}
\frac{1}{\sqrt{a_{00}}}=1/\sqrt{-\frac{\cos^2{\frac{\pi}{w}}\cos^2{\frac{\pi}{u}}-\cos^2{\frac{\pi}{w}}+\sin^2{\frac{\pi}{v}}-\cos^2{\frac{\pi}{u}}}{\cos^2{\frac{\pi}{w}}-\sin^2{\frac{\pi}{v}}}}, \\
\frac{1}{\sqrt{a_{33}}}=1/\sqrt{-\frac{\cos^2{\frac{\pi}{w}}
\cos^2{\frac{\pi}{u}}-\cos^2{\frac{\pi}{w}}+\sin^2{\frac{\pi}{v}}
-\cos^2{\frac{\pi}{u}}}{\cos^2{\frac{\pi}{u}}-\sin^2{\frac{\pi}{v}}}}. 
\end{gathered} \notag
\end{equation}
We obtain the statement of the lemma directly from the above formulas. \ \ \ $\square$
\medbreak
We may assume that $u \le w$ because of the symmetrical role of parameters $u$ and $w$. Therefore,
$d(C,A_1) \le d(Q,A_2)$ inequality holds.

We have to distinguish $2$ different main cases 
all of which we set up from the optimal congruent hyperball packing configuration 
described in the former subsection.
\begin{enumerate}
\item $h=d(H,J)/2 \le d(C,A_1)$ (and so $h \le d(Q,A_2)$. 
\smallbreak
Both polar planes are assigned hyperspheres.
It is clear, that in the ''starting configuration" (congruent case) the heights 
of optimal hyperballs $\mathcal{H}^{h^{i}}(\mathcal{A}_i)$ $(i=0,3)$ 
are equal $h=h^{0}=h^3$.
\begin{enumerate}
\item We consider the hyperball $\mathcal{H}^{h^{0}}(\mathcal{A}_0)$ 
and blow up this hyperball (hypersphere) keeping the hyperballs 
$\mathcal{H}^{h^3}(\mathcal{A}_3)$ tangent to it upto this hypersphere 
touch the plane $A_1A_2A_3$ or the plane $EJQ$ (see Fig.~1). 
During this expansion the height of hyperball 
$\mathcal{H}^{h^0}$ is  
$h^0=h+x$ where $x \in [0,\min\{d(C,A_1)-h, h \}]$. 
The height of hyperball $\mathcal{H}^{h^3}(\mathcal{A}_3)$ is $h^3=h-x$. 
(If $x=0$ then the hyperballs are congruent.) 
\item We consider the hyperball $\mathcal{H}^{h^{3}}(\mathcal{A}_3)$ 
and blow up this hyperball (hypersphere) keeping the hyperballs 
$\mathcal{H}^{h^0}(\mathcal{A}_0)$ tangent to it upto this hypersphere 
touch the plane $A_0A_1A_2$ or the plane $CLH$ (see Fig.~1). 
During this expansion the height of hyperball
$\mathcal{H}^{h^3}(\mathcal{A}_3)$ is  
$h^3=h+x$ where $x \in [0,\min\{d(Q,A_2)-h, h \}]$. 
The height of hyperball $\mathcal{H}^{h^0}(\mathcal{A}_0)$ is $h^0=h-x$.  
\end{enumerate}
\item $h=d(C,A_1)$, $(d(C,A_1) \le d(Q,A_2)$ because $u \le w$.
\smallbreak
In this case we distinguish two subcases:
\begin{enumerate}
\item We blow up the hyperball $\mathcal{H}^{h^{3}}(\mathcal{A}_3)$ 
upto this hypersphere touch the plane $A_0A_1A_2$ or the hypersphere 
$\mathcal{H}^{h^{0}}(\mathcal{A}_0)$ (see Fig.~1). 
During this expansion the height of hyperball 
$\mathcal{H}^{h^0}(\mathcal{A}_0)$ is  
$h^0=h$ constant and the height of hyperball $\mathcal{H}^{h^3}(\mathcal{A}_3)$ 
is $h^3=h+x$ where $x \in [0,\min\{d(Q,A_2)-d(C,A_1), d(H,J)-2h \}]$.
(If $x=0$ then the hyperballs are congruent.) 
\item If in the above situation the hyperball $\mathcal{H}^{h^3}(\mathcal{A}_3)$
first touches the hyperball $\mathcal{H}^{h^0}(\mathcal{A}_0)$ 
$(h^3=d(H,J)-h \le d(Q,A_2))$ then we can continue the
blowing of hyperball $\mathcal{H}^{h^3}(\mathcal{A}_3)$ upto 
this hypersphere 
touch the plane $A_0A_1A_2$ or the plane $CLH$ (see Fig.~1). 
During this expansion the height of hyperball
$\mathcal{H}^{h^3}(\mathcal{A}_3)$ is  
$h^3=d(H,J)-h+x$ where $x \in [0,\min\{d(Q,A_2)-d(H,J)+h, h \}]$. 
The height of hyperball $\mathcal{H}^{h^0}(\mathcal{A}_0)$ is $h^0=h-x$.
\end{enumerate}
\end{enumerate}
We extend this arrangements to images of the hyperballs 
$\mathcal{H}^{h^i}(\mathcal{A}_i)$ $(i\in\{0,3\})$ by the considered 
Coxeter group and obtain non-congruent hyperball packing $\mathcal{B}(x)$. 
Its density is defined by the following
\begin{definition} The locally density functions $\delta^j_x(\mathcal{O}(u,v,w)$ 
of the non-cong\-ru\-ent hyperball packings related to 
$\mathcal{O}(u,v,w)$ in the above cases ($j\in \{1,2\}$) are 
defined by next formulas:
\begin{equation}
\begin{gathered}
\delta^j_x(\mathcal{O}(u,v,w)):=\\ = \frac{Vol(\mathcal{H}^{h^0}(\mathcal{A}_0)\cap \mathcal{O}(u,v,w))+
Vol(\mathcal{H}^{h^3}(\mathcal{A}_3)\cap \mathcal{O}(u,v,w))}{Vol({\mathcal{O}(u,v,w)})}, \notag
\end{gathered}
\end{equation}
where $0 \le h^0,h^3$ are suitable real parameters related to the above main non-congruent cases (depend on $x$ parameter) 
and $\mathcal{A}_0$ is the area of the triangle $CLH$ and 
$\mathcal{A}_3$ is the area of the triangle $JEQ$ (see (4.1)).
\end{definition}
{\it {\bf The main problem is}: What is the maximum of density function 
$\delta^j_x(\mathcal{O}(u,$ $v,w))$
for suitable integer parameters $u,v,w$ where
$x \in \mathbb{R}$, and $x \in [0,\min\{$ $d(Q,A_2)-d(C,A_1), d(H,J)-2h \}]$ 
or $x \in [0,\min\{d(Q,A_2)-d(H,J)+h, h \}]$} (see the above two main cases).
\subsubsection{Numerical data of non-congruent hyperball packing arrangements}
First, we illustrate our computation method for given parameters and 
then we summarize the numerical data of optimal non-congruent hyperball arrangements 
for several given parameters in Table 3.
\medbreak
{\bf Results for parameters $u=5$, $v=4$, $w=5$}:
\begin{figure}[htb]
\centering
\includegraphics[width=50mm]{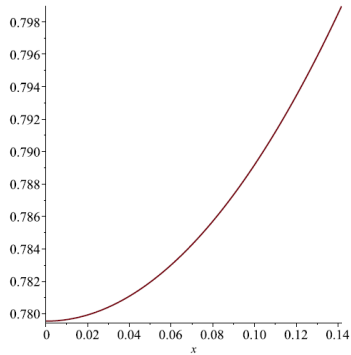} 
\caption{Density function $\delta^1_x(\mathcal{O}(5,4,5))$ where
$x \in [0, d(C,A_1)-h \approx 0.14166]$}
\label{}
\end{figure}
\begin{equation}
\begin{gathered}
d(A_1,C)=\mathrm{arcosh}{\frac{1}{\sqrt{a_{00}}}}= 
d(A_2,Q)=\mathrm{arcosh}{\frac{1}{\sqrt{a_{33}}}} \approx 1.02221,\\
d(J,H)/2=\frac{1}{2}\mathrm{arcosh}{\frac{-a_{03}}{\sqrt{a_{00}a_{33}}}}\approx 0.88055
\end{gathered} \notag 
\end{equation}
Therefore, the heights in the starting position ($x=0$): $h=h^0=h^3=d(J,H)/2$ $ \approx 0.88055$ 
(see (4.1), (4.2), (4.3)) and $Vol({\mathcal{O}(5,4,5)}) \approx 0.38325$.

We consider the hyperball $\mathcal{H}^{h^{0}}(\mathcal{A}_0)$ 
and blow up this hyperball (hypersphere) keeping the hyperballs 
$\mathcal{H}^{h^3}(\mathcal{A}_3)$ tangent to it upto this hypersphere 
touch the plane $A_1A_2A_3$ (see Fig.~1). 
During this expansion the height of hyperball 
$\mathcal{H}^{h^0}$ is  
$h^0=h+x$ where . 
The height of hyperball $\mathcal{H}^{h^3}(\mathcal{A}_3)$ is $h^3=h-x$. 
\begin{figure}[htb]
\centering
\includegraphics[width=110mm]{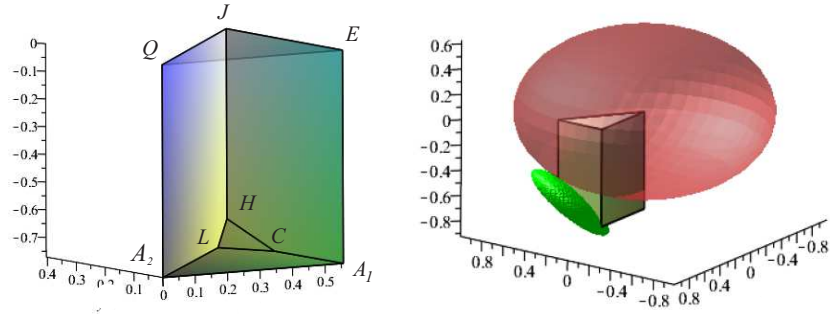} 
\caption{Locally optimal non-congruent hyperball packing configuration related to parameters $\{5,4,5\} $}
\label{}
\end{figure}
We get by the Definition 4.7 of the density function $\delta^1_x(\mathcal{O}(5,4,5))$ that it is a strictly increasing function in the intervall 
$[0, d(C,A_1)-h \approx 0.14166]$ (see Fig.~4). Thus, the optimal arrangement belongs to the parameter $x^{opt}=d(C,A_1)-h \approx 0.14166$:
\begin{equation} 
\begin{gathered}
\delta^1_{x^{opt}}(\mathcal{O}(5,4,5))=\\ 
= \frac{Vol(\mathcal{H}^{h^0}(\mathcal{A}_0)\cap \mathcal{O}(5,4,5))+
Vol(\mathcal{H}^{h^3}(\mathcal{A}_3)\cap \mathcal{O}(5,4,5))}
{Vol({\mathcal{O}(5,4,5)})} \approx 0.79895, \notag
\end{gathered}
\end{equation}
where $x^{opt}=d(C,A_1)-h \approx 0.14166$ and 
$\mathcal{A}_0=\mathcal{A}_3=\pi/2-\pi/5-\pi/4$. (see (4.1)).
\medbreak

\medbreak
{\scriptsize
\centerline{\vbox{
\halign{\strut\vrule~\hfil $#$ \hfil~\vrule
&\quad \hfil $#$ \hfil~\vrule
&\quad \hfil $#$ \hfil\quad\vrule
&\quad \hfil $#$ \hfil\quad\vrule
&\quad \hfil $#$ \hfil\quad\vrule
\cr
\noalign{\hrule}
\multispan5{\strut\vrule\hfill\bf Table 3, allowing non-congruent hyperballs  \hfill\vrule}%
\cr
\noalign{\hrule}
\noalign{\vskip2pt}
\noalign{\hrule}
\{u,v,w\} & h^0  & h^3 & \sum_{i=0,3} {Vol}(\mathcal{H}^{h^i}(\mathcal{A}_i))\\
& \delta^1({\mathcal{O}(u,v,w)})\cr
\noalign{\hrule}
\{7,3,7\} & 1.23469 & 1.23469 & 0.31172 & 0.81335 \cr
\noalign{\hrule}
\{7,3,8\} & 0.93100 & 1.25596 & 0.32520 & 0.78690 \cr
\noalign{\hrule}
\{7,3,9\} & 0.76734 & 1.27042 & 0.32892 & 0.76189 \cr
\noalign{\hrule}
\vdots & \vdots  & \vdots  & \vdots  & \vdots  \cr
\noalign{\hrule}
\{7,3,50\} & 0.11380 & 1.32226 & 0.23307 & 0.47549 \cr
\noalign{\hrule}
\vdots & \vdots  & \vdots  & \vdots  & \vdots  \cr
\noalign{\hrule}
\{8,3,8\} & 0.94946 & 0.94946 & 0.33794 & 0.76143 \cr
\noalign{\hrule}
\{8,3,9\} & 0.78366 & 0.96206 & 0.34107 & 0.73718 \cr
\noalign{\hrule}
\{8,3,10\} & 0.67409 & 0.97104 & 0.33990 & 0.71504 \cr
\noalign{\hrule}
\vdots & \vdots  & \vdots  & \vdots  & \vdots  \cr
\noalign{\hrule}
\{8,3,50\} & 0.11668 & 1.00753 & 0.24051 & 0.46032 \cr
\noalign{\hrule}
\vdots & \vdots  & \vdots  & \vdots  & \vdots  \cr
\noalign{\hrule}
\{5,4,5\} & 0.73890 & 1.02221 & 0.36903 & 0.79895 \cr
\noalign{\hrule}
\{5,4,6\} & 0.73969 & 0.83611 & 0.39956 & 0.78736 \cr
\noalign{\hrule}
\{5,4,7\} & 0.59326 & 0.90486 & 0.41263 & 0.77517 \cr
\noalign{\hrule}
\vdots & \vdots  & \vdots  & \vdots  & \vdots  \cr
\noalign{\hrule}
\{5,4,50\} & 0.07206 & 1.19054 & 0.35623 & 0.60082 \cr
\noalign{\hrule}
\vdots & \vdots  & \vdots  & \vdots  & \vdots  \cr
\noalign{\hrule}
\{4,5,4\} & 0.55565 & 1.06128 & 0.34184 & 0.79382 \cr
\noalign{\hrule}
\{4,5,5\} & 0.69129 & 0.69129 & 0.38284 & 0.76893 \cr
\noalign{\hrule}
\{4,5,6\} & 0.53064 & 0.77568 & 0.39374 & 0.74331 \cr
\noalign{\hrule}
\vdots & \vdots  & \vdots  & \vdots  & \vdots  \cr
\noalign{\hrule}
\{4,5,50\} & 0.05502 & 1.13842 & 0.32720 & 0.55161 \cr
\noalign{\hrule}
\{5,5,5\} & 0.35764 & 0.77537 & 0.41589 & 0.72618 \cr
\noalign{\hrule}}}}}
Similarly to the above cases the volume of a doubly truncated orthoscheme can be computed by Theorem 2.1 and here we allow non-congruent hyperballs. 
The volume sum $\sum_{i=0,3} {Vol}(\mathcal{H}^{h^i}(\mathcal{A}_i))$ of hyperball 
pieces lying in ${\mathcal{O}(u,v,w)}$ can be computed by the formula (2.3), (4.1), (4.2), (4.3) and by the above described computation method 
for each given possible parameters $u,v,w$. But, the computations can be contained some subcases so, the determination of the
densest hyperball configuration for given parameters $u,v,w$ more complicated as the above example shows. 

The maximal density of the congruent hyperball packing related to the doubly truncated orthoscheme tilings with ''non-congruent hyperballs" 
-- $\delta^i_x({\mathcal{O}(u,v,w)})$ (see Definition 4.7, $i\in\{1,2\}$) -- can be computed for each possible parameters.

Finally, we can formulate only the next Theorem and Conjecture (see Fig.~5):
\begin{theorem}
The density functions $\delta^i_x({\mathcal{O}(u,v,w)})$, ($\frac{1}{u}+\frac{1}{v} < \frac{1}{2},~\frac{1}{v}+\frac{1}{w} < \frac{1}{2}$. $3 \le u,v,w \in \mathbb{N}$
and $x \in [0,\min\{d(Q,A_2)-d(C,A_1), d(H,J)-2h \}]$ 
or $x \in [0,\min\{d(Q,A_2)-d(H,J)+h, h \}]$) attain their maximum at $\{u,v,w\}=\{7,3,7\}$ with $x=0$ (congruent case) among the investigated parameters with 
density $\delta^1_x({\mathcal{O}(7,3,7)}) \approx 0.81335$ (see Table 3). 
\end{theorem}
\begin{conjecture}
The density functions $\delta^i_x({\mathcal{O}(u,v,w)})$, ($\frac{1}{u}+\frac{1}{v} < \frac{1}{2},~\frac{1}{v}+\frac{1}{w} < \frac{1}{2}$. $3 \le u,v,w \in \mathbb{N}$
and $x \in [0,\min\{d(Q,A_2)-d(C,A_1), d(H,J)-2h \}]$ 
or $x \in [0,\min\{d(Q,A_2)-d(H,J)+h, h \}]$) attain their maximum at $\{u,v,w\}=\{7,3,7\}$ with $x=0$ (congruent case) with 
density $\delta^1({\mathcal{O}(7,3,7)}) \approx 0.81335$ (see Table 3). 
\end{conjecture}
\subsubsection{On non-extendable congruent hyperball packings to parameters $\{ p,3,p\}$ $(6 < p < 7, ~ p\in \mathbb{R})$}
We can construct infinitely many congruent and non congruent hyperball configuration whose densities are locally larger, that the 
B\"or\"oczky-Florian density upper bound $(\approx 0.85328)$. Now, we describe only one, congruent locally dense hyperball arrangement related to parameters 
$\{ p,3,p\}$ $(6 < p < 7, ~ p\in \mathbb{R})$.

The computation method described in the former sections is suitable to determine the densities of congruent hyperball packings for
parameters $(6 < p < 7, ~ p \in \mathbb{R})$ as well. To each $p$ parameter belongs a doubly truncated orthoscheme and therefore we can determine similarly to the
above cases the corresponding maximal density of its optimal congruent hyperball packing. But these packings can not be extended to the 3-dimensional
space. Analysing these non-extendable packings for parameters $(6 < p < 7, ~ p \in \mathbb{R})$ we obtain the following (see Fig.~6)
\begin{theorem}
The function $\delta^1({\mathcal{O}(u=p,v=3,w=p)})$, ($\frac{1}{u}+\frac{1}{v} < \frac{1}{2},~\frac{1}{v}+\frac{1}{w} < \frac{1}{2}$. $3 \le u,v,w \in \mathbb{N}$, see Definition 4.1)
is attained its maximum for the parameter $p_{opt} \approx 6.05061$
and the density is  $\approx 0.85461$. That means that these
locally optimal hyperball configuration provides larger densities that the B\"or\"oczky-Florian density upper bound $(\approx 0.85328)$ for ball and
horoball packings (\cite{B--F64}).
\end{theorem}
\begin{figure}[htb]
\centering
\includegraphics[width=50mm]{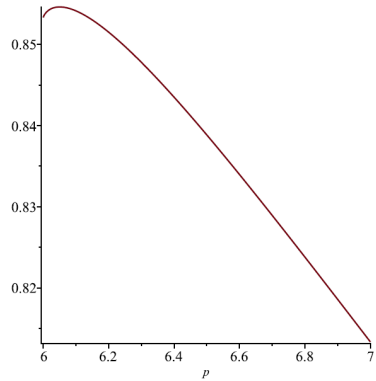} 
\caption{Locally optimal non-extendable congruent hyperball packing configuration related to parameters $\{ p,3,p\}$ $(6 < p < 7, ~ p\in \mathbb{R})$}
\label{}
\end{figure}
In hyperbolic spaces $\HYN$ ($n \ge 3$) the  problems of the densest horoball and hyperball packings have not been settled yet, 
in general (see e.g. \cite{KSz14}, \cite{Sz12}, \cite{Sz12-2}).
Moreover, the optimal sphere packing problem can be extended to the other homogeneous Thurston geometries, e.g. $\NIL$, $\SOL$, $\SLR$. 
For these non-Euclidean geometries only very few results are known (e.g. \cite{Sz14-1} and the references given there).

By the above these we can say that the revisited Kepler problem keep several interesting open questions.
Detailed studies are the objective of ongoing research. Applications of the above projective method seem to be interesting in (non-Euclidean)
crystallography as well.


\noindent
\footnotesize{Budapest University of Technology and Economics, Institute of Mathematics, \\
Department of Geometry, \\
H-1521 Budapest, Hungary. \\
E-mail:~szirmai@math.bme.hu \\
http://www.math.bme.hu/ $^\sim$szirmai}


\begin{thebibliography}{999999999}

\bibitem{Be} Bezdek,~K.
Sphere Packings Revisited,
\textit{Eur. J. Combin.}, {\bf{27/6}} (2006), \rm 864--883.
%
\bibitem{B31} Bolyai,~J.
\textit{Appendix, Scientiam spatii absolute veram exhibens},
\rm Marosvásárhely, (1831).
%
\bibitem{B78} B\"or\"oczky,~K.
Packing of spheres in spaces of constant curvature,
\textit{Acta Math. Acad. Sci. Hungar.}, {\bf{32}} (1978), \rm 243--261.
%
\bibitem{B--F64} B\"or\"oczky,~K. ~---~ Florian,~A.
\"Uber die dichteste Kugelpackung im hyperbolischen Raum, \textit{Acta Math. Acad. Sci. Hungar.},
{\bf{15}} (1964), \rm 237--245.
%
\bibitem{G--K--K} Fejes~T\'oth,~G.~---~Kuperberg,~G.~---~Kuperberg,~W.
Highly Saturated Packings and Reduced Coverings,
\textit{Monatsh. Math.}, {\bf{125/2}} (1998), \rm 127--145.
%
\bibitem{FTL} Fejes~T\'oth,~L.
Regular Figures,
\textit{Macmillan (New York)}, 1964.
%
\bibitem{IH90} Im Hof,~H.-C. Napier cycles and hyperbolic Coxeter groups,
\textit{Bull. Soc. Math. Belgique}, {\bf{42}} (1990), 523--545.
%
\bibitem{K89} Kellerhals,~R. On the volume of hyperbolic polyhedra,
\textit{Math. Ann.}, {\bf{285}} (1989), 541--569.
%
\bibitem{KSz14} Kozma,~R.T.~---~Szirmai,~J.
New Lower Bound for the Optimal Ball Packing Density of Hyperbolic 4-space,
\emph{Discrete Comput. Geom.}, {\bf 53/1} (2015), 182-198, DOI: 10.1007/s00454-014-9634-1.
%
\bibitem{Mol97} Moln\'ar,~E.
The Projective Interpretation of the eight 3-dimensional homogeneous geometries,
\textit{Beitr. Algebra Geom.,} {\bf{38/2}} (1997), \rm 261--288.
%
\bibitem{MoSzi18} Moln\'ar,~E.~---~Szirmai,~J.
Top dense hyperbolic ball packings and coverings for complete Coxeter orthoscheme groups,
\textit{Publications de l'Institut Mathématique}, {\bf 103(117)}  (2018), \rm 129--146, 
DOI: 10.2298/PIM1817129M,  arXiv: 1612.04541.
%
\bibitem{S14} Stojanovi\' c,~M.
Coxeter Groups as Automorphism Groups
of Solid Transitive 3-simplex Tilings,
\textit{Filomat,} {\bf{28/3}} (2014), 557--577, DOI 10.2298/FIL1403557S.
%
\bibitem{S17} Stojanovi\' c,~M.
Hyperbolic space groups and their supergroups for fundamental simplex tilings,
\textit{Acta Math. Hungar.,} {\bf{153/2}} (2017), 276--288, DOI: 10.1007/s10474-017-0761-z.
%
\bibitem{Sz14} Szirmai,~J. Hyperball packings in hyperbolic $3$-space,
\textit{Mat. Vesn.}, {\bf 70/3} (2018), 211--221.
%
\bibitem{Sz17} Szirmai,~J. Packings with horo- and hyperballs generated by simple frustum orthoschemes,
\textit{Acta Math. Hungar.}, {\bf 152/2} (2017), 365--382, DOI:10.1007/s10474-017-0728-0.
%
\bibitem{Sz17-1} Szirmai,~J. Density upper bound of congruent and non-congruent hyperball packings generated by truncated regular simplex tilings,
\textit{Rendiconti del Circolo Matematico di Palermo Series 2}, {\bf 67} (2018),  307--322,
DOI: 10.1007/s12215-017-0316-8,  arXiv:1510.03208.
%
\bibitem{Sz17-2} Szirmai,~J. Decomposition method related to saturated hyperball packings,
\textit{Submitted manuscript}, (2017).
%
\bibitem{Sz07-1} Szirmai,~J.
The optimal ball and horoball packings to the Coxeter honeycombs  in the hyperbolic $d$-space,
\textit{Beitr. Algebra Geom.,} {\bf{48/1}} (2007), 35--47.
%
\bibitem{Sz12} Szirmai,~J.
Horoball packings to the totally asymptotic regular simplex in the hyperbolic $n$-space,
\emph{Aequat. Math.},  {\bf 85} (2013), 471-482,  DOI: 10.1007/s00010-012-0158-6.
%
\bibitem{Sz12-2}Szirmai,~J.
Horoball packings and their densities by generalized simplicial density function in the hyperbolic space,
\emph{Acta Math. Hungar.,}
{\bf 136/1-2} (2012), 39--55, DOI: 10.1007/s10474-012-0205-8.
%
\bibitem{Sz06-1} Szirmai,~J. The $p$-gonal prism tilings and their optimal hypersphere packings in the hyperbolic
3-space,
\textit{Acta Math. Hungar.}, {\bf{111 (1-2)}} (2006), 65--76.
%
\bibitem{Sz06-2} Szirmai,~J. The regular prism tilings and their optimal hyperball packings in the hyperbolic $n$-space,
\textit{Publ. Math. Debrecen}, {\bf{69 (1-2)}} (2006), 195--207.
%
\bibitem{Sz13-3} Szirmai,~J. The optimal hyperball packings related to the smallest compact arithmetic $5$-orbifolds,
\textit{Kragujevac J. Math.} {\bf 40(2)} (2016), 260-270, DOI:10.5937/KgJMath1602260S.
%
\bibitem{Sz13-4} Szirmai,~J. The least dense hyperball covering to the regular prism tilings in the hyperbolic $n$-space,
\textit{Ann. Mat. Pur. Appl.} {\bf 195/1} (2016), 235--248, DOI: 10.1007/s10231-014-0460-0.
%
\bibitem{Sz14-1}
{Szirmai,~J.}
A candidate for the densest packing with equal balls in Thurston geometries,
\emph{Beitr. Algebra Geom.}, {\bf 55/2} (2014), 441- 452, DOI: 10.1007/s13366-013-0158-2.
%
\bibitem{Sz18} Szirmai,~J. Hyperball packings related to octahedron and cube tilings in hyperbolic space,
\textit{Submitted manuscript}, (2018).
%
\bibitem{V72} Vermes,~I. \"Uber die Parkettierungsm\"oglichkeit des dreidimensionalen hyperbolischen
Raumes durch kongruente Polyeder,
\textit{Studia Sci. Math. Hungar.}, {\bf{7}} (1972), 267--278.
%
\bibitem{V79} Vermes,~I. Ausf\"ullungen der hyperbolischen Ebene durch kongruente Hyperzykelbereiche,
\textit{Period. Math. Hungar.}, {\bf{10/4}} (1979), 217--229.
%
\bibitem{V81} Vermes,~I. Über regul\"are Überdeckungen der Bolyai-Lobatschewskischen Ebene durch kongruente Hyperzykelbereiche,
\textit{Period. Math. Hungar.}, {\bf{25/3}} (1981), 249--261.
\end{thebibliography}
\end{document}